\documentclass{article}
\usepackage{amsmath,amssymb}
\newcommand{\Z}{\ensuremath{\mathbf Z}}
\newcommand{\N}{\ensuremath{ \mathbf N }}
\newtheorem{theorem}{Theorem}
\newcommand{\bt}{\begin{theorem}}
\newcommand{\et}{\end{theorem}}
\newcommand{\pf}{{\bf Proof}.\ }
\newcommand{\beq}{\begin{equation}}
\newcommand{\eeq}{\end{equation}}
\newcommand{\benum}{\begin{enumerate}}
\newcommand{\eenum}{\end{enumerate}}
\DeclareMathOperator{\card}{card}
\newcommand{\eop}{$\square$\vspace{.3cm}}

\begin{document}
\title{Generalized additive bases, K\" onig's lemma,\\
and the Erd\H os-Tur\' an conjecture\footnote{2000 Mathematics
Subject Classification:  11B13, 11B34, 11B05.
Key words and phrases.  Additive bases, sumsets, representation functions, 
Erd\H os-Tur\' an conjecture, K\" onig's lemma.}}
\author{Melvyn B. Nathanson\thanks{This work was supported
in part by grants from the NSA Mathematical Sciences Program
and the PSC-CUNY Research Award Program.}\\
Department of Mathematics\\
Lehman College (CUNY)\\
Bronx, New York 10468\\
Email: nathansn@alpha.lehman.cuny.edu}
\maketitle

\begin{abstract}
Let $A$ be a set of nonnegative integers.
For every nonnegative integer $n$ and positive integer $h$,
let $r_{A}(n,h)$ denote
the number of representations of $n$ in the form
$ n =  a_1 + a_2 + \cdots +  a_h,$ where $a_1, a_2, \ldots,a_h \in A$ 
and $a_1 \leq a_2 \leq \cdots \leq a_h.$
The infinite set $A$ is called a {\em basis of order $h$} if 
$r_{A}(n,h) \geq 1$ for every nonnegative integer $n$. 
Erd\H os and Tur\' an conjectured that 
$\limsup_{n\rightarrow\infty}r_A(n,2) = \infty$
for every basis $A$ of order $2$.
This paper introduces a new class of additive bases
and a general additive problem, a special case of which is the
Erd\H os-Tur\' an conjecture.
K\" onig's lemma on the existence of infinite paths in certain graphs 
is used to prove that this general problem is equivalent to a related
problem about finite sets of nonnegative integers.
\end{abstract}

\section{Representation functions and the Erd\H os-Tur\' an conjecture}
Let $\N_0$ and \Z\ denote the nonnegative integers and integers, 
respectively.
Let $A$ be a finite set of integers.  
We denote the largest element of $A$ by $\max(A)$ 
and the cardinality of $A$ by $\card(A)$.
For any real numbers $a$ and $b$, we denote by $[a,b]$ 
the finite set of integers $n$ such that $a \leq n \leq b.$

For any set $A$ of integers, we denote by $r_A(n,h)$ the number
of representations of $n$ in the form $n = a_1 + a_2 + \cdots + a_h,$ 
where $a_1, a_2, \ldots,a_h \in A$ and $a_1 \leq a_2 \leq \cdots \leq a_h.$
The function $r_{A}$ is called the {\em unordered representation function} 
of the set $A$, or, simply, the {\em representation function} 
of $A$.

The set $A$ of nonnegative integers is called 
a {\em basis of order $h$}
if every nonnegative integer can be represented as the sum of $h$
not necessarily distinct elements of $A$.
If $A$ is a basis of order $h$ with representation function $r_A,$
then 
\[
1 \leq r_A(n,h) < \infty
\]
for all nonnegative integers $n$.
We call $A$ an {\em asymptotic basis of order $h$} if $r_A(n,h) = 0$
for only finitely many $n \in \N_0.$
If $h \geq 2$ and $f: \N_0 \rightarrow \N_0$ is any function 
such that $\card(f^{-1}(0)) < \infty$, then
the {\em representation function problem} is to determine
if there exists a set $A$ of nonnegative integers such that $r_A(n,h) = f(n)$
for all $n \geq 0$.

In the case of additive bases for the set of all integers, 
Nathanson~\cite{nath03c} proved that every function is a representation function, 
that is, if $f:\Z \rightarrow \N_0 \cup \{\infty\}$ 
satisfies the condition $\card(f^{-1}(0)) < \infty$, 
then for every $h \geq 2$ there exists a set $A$ of integers 
such that $r_A(n,h) = f(n)$ for every integer $n$. 

A special case of the representation function problem for nonnegative integers with $h=2$ 
is the conjecture of Erd\H os and Tur\' an~\cite{erdo-tura41} that 
the representation function $r_A(n,2)$ of an asymptotic basis $A$ 
of order 2 must be unbounded, that is,
\[
\liminf_{n\rightarrow\infty} r_A(n,2) > 0 
\Longrightarrow \limsup_{n\rightarrow\infty} r_A(n,2) = \infty.
\]
This is an important unsolved problem in additive number theory.

Dowd~\cite{dowd88} and Grekos, Haddad, Helou, and Pikho~\cite{grek-hadd-helo-pihk03}
have given various equivalent formulations of the Erd\H os-Tur\' an conjecture.
In particular, Dowd proved that there exists a set $A$ of nonnegative integers
and a number $c$ such that $r_A(n,2) \in [1,c]$ for all nonnegative integers $n$
if and only if for every $N$ there exists a finite set $A_N$ 
of nonnegative integers with $\max(A_N) \geq N$ and
$r_A(n,2) \in [1,c]$  for all $n = 0,1,\ldots,\max(A_N).$
In this paper we apply Dowd's method to obtain similar results for 
a new class of generalized additive bases.

\section{Generalized additive bases}
We extend the idea of an additive basis of order $h$ as follows:
Let $\mathcal{H} = \{H_n\}_{n=0}^{\infty}$ 
be a sequence of nonempty finite sets of positive integers.
For any set $A$ of nonnegative integers, we define
the representation function
\[
r_A(n,H_n) = \sum_{h_n \in H_n} r_{A}(n, h_n).
\]
The set $A$ of nonnegative integers will be called a 
{\em basis of order $\mathcal{H}$} if 
\beq  \label{konig:rinq}
r_{A}(n, H_n) \geq 1
\eeq
for all $n \geq 0,$ and an {\em asymptotic basis of order $\mathcal{H}$} 
if the representation function satisfies~{(\ref{konig:rinq}) 
for all sufficiently large $n$.

Let $\mathcal{R} = \{R_n\}_{n=0}^{\infty}$ 
be a sequence of nonempty finite sets of positive integers.
If 
\beq       \label{konig:rineq2}
r_A(n,H_n) \in R_n 
\eeq
for every nonnegative integer $n$, then $A$ will be called an 
{\em $\mathcal{R}$-basis of order $\mathcal{H}$}.
Since each $R_n$ is a nonempty set of positive integers,
it follows that every $\mathcal{R}$-basis of order $\mathcal{H}$
is a basis of order $\mathcal{H}$.
We shall call the set $A$ an {\em asymptotic $\mathcal{R}$-basis of order $\mathcal{H}$}
if~(\ref{konig:rineq2}) holds for all sufficiently large $n$.

For any sequences $\mathcal{H} = \{H_n\}_{n=0}^{\infty}$ 
and $\mathcal{R} = \{R_n\}_{n=0}^{\infty}$ 
of nonempty finite sets of positive integers,
we can ask if there exists an $\mathcal{R}$-basis of order $\mathcal{H}$
or an asymptotic $\mathcal{R}$-basis of order $\mathcal{H}$.
This is the {\em generalized representation function problem.}
The original Erd\H os-Tur\' an conjecture corresponds to the special case $H_n = \{ 2\}$ 
and $R_n = [1,c]$ for all $n \geq 0$.
It is an open problem to determine the number
of distinct $\mathcal{R}$-bases of order $\mathcal{H}$  
and asymptotic $\mathcal{R}$-bases of order $\mathcal{H}$ 
for a given pair of sequences $\mathcal{H}$ and $\mathcal{R}.$ 

Let $h \geq 2$ and let $f$ be a function such that $f(n)$ is a positive 
integer for every nonnegative integer $n.$
We introduce the sets $H_n = \{h\}$ and $R_n = \{f(n)\}$ for all $n$, 
and the sequences $\mathcal{H} = \{H_n\}_{n=0}^{\infty}$ 
and $\mathcal{R} = \{R_n\}_{n=0}^{\infty}$.  
Then an $\mathcal{R}$-basis of order $\mathcal{H}$ 
is a basis $A$ of order $h$
whose representation function satisfies $r_A(n,h) = f(n)$
for all $n \in \N_0$, and so the 
representation function problem for bases of order $h$ is  a special case
of the generalized representation function problem.

An $\mathcal{R}$-basis of order $\mathcal{H}$
is not necessarily infinite.  For example, 
if $H_0 = R_0 = \{1\}$ and if $H_n = \{n\}$ and $1 \in R_n$
for all $n \geq 1,$ then the set $\{0,1\}$ is an
$\mathcal{R}$-basis of order $\mathcal{H}$.

\bt          \label{konig:theorem:finite}
Let $\mathcal{H} = \{H_n\}_{n=0}^{\infty}$ 
be a sequence of nonempty finite sets of positive integers.
There exists a finite set $A$ that is a basis of order $\mathcal{H}$
or an asymptotic basis of order $\mathcal{H}$ if and only if
\[
\liminf_{n\rightarrow\infty} \frac{\max(H_n)}{n} > 0.
\]
\et

\pf
Let $h_n^* = \max(H_n).$
Let $A$ be a finite set of nonnegative integers that is 
a basis of order $\mathcal{H}$.  Then $0,1 \in A$ and so $\max(A)\geq 1.$
Every positive integer $n$ can be represented as the sum 
of $h_n$ elements of $A$ for some $h_n \in H_n,$ 
and so 
\[
n \leq h_n\max(A)\leq h_n^*\max(A).
\]
It follows that 
\[
\liminf_{n\rightarrow\infty} \frac{h_n^*}{n} \geq \frac{1}{\max(A)} > 0.
\]

Conversely, if $\liminf_{n\rightarrow\infty} h_n^*/n > 0,$
then there exists a positive integer $m$ such that 
\[
h_n^* \geq \frac{n}{m}
\]
for all $n \geq 0$.
Consider the finite set $A = [0,m].$  By the division algorithm,
every positive integer $n$ can be written in the form
$n = qm+r,$ where $q$ and $r$ are nonnegative integers 
and $0 \leq r \leq m-1.$
If $r=0,$ then $q = n/m \leq h_n^*$ and 
\[
n = q\cdot m + (h_n^*-q)\cdot 0 \in h_n^*A.
\]
If $1 \leq r \leq m-1,$ then $q = (n-r)/m < h_n^*$.  
Since $h_n^*$ and $q$ are integers,
it follows that $h_n^* \geq q+1$ and 
\[
n = q\cdot m + 1\cdot r + (h_n^*-q-1)\cdot 0 \in h_n^*A.
\]
In both cases, $r_A(n,H_n) \geq r_A(n,h_n^*) \geq 1$, and
the finite set $A$ is a basis of order $\mathcal{H}$.

If $A$ is an asymptotic basis of order $\mathcal{H},$
then $r_A(n,H_n) = 0$ for only finitely many $n \in \N_0,$
and so there is a finite set $F$ of nonnegative integers
such that $A\cup F$ is a basis of order $\mathcal{H}$.
Therefore, there exists a finite set that is a basis 
of order $\mathcal{H}$ if and only if there exists a finite set that is 
an asymptotic basis of order $\mathcal{H}$.
This completes the proof.
\eop

Let $\mathcal{H} = \{H_n\}_{n=0}^{\infty}$ 
and $\mathcal{R} = \{R_n\}_{n=0}^{\infty}$ 
be sequences of nonempty finite sets of positive integers.
A nonempty finite set $A$ of nonnegative integers 
will be called a {\em finite basis of order $\mathcal{H}$} 
if 
\[
r_A(n,H_n) \geq 1 
\]
for all $n \in \left[ 0,\max(A)\right]$,
and a {\em finite $\mathcal{R}$-basis of order $\mathcal{H}$} 
if 
\[
r_A(n,H_n) \in R_n 
\]
for all $n \in \left[0,\max(A)\right]$.

\bt      \label{konig:theorem:facts}
Let $\mathcal{H} = \{H_n\}_{n=0}^{\infty}$ 
and $\mathcal{R} = \{R_n\}_{n=0}^{\infty}$ 
be sequences of nonempty finite sets of positive integers.
\benum
\item[(i)]
If $A$ is a basis of order $\mathcal{H}$,
or if $A$ is a finite basis of order $\mathcal{H}$
with $\max(A) \geq 1,$ then $0,1 \in A.$
\item[(ii)]
If $A$ is an $\mathcal{R}$-basis of order $\mathcal{H}$
or if $A$ is a finite $\mathcal{R}$-basis of order $\mathcal{H}$
with $\max(A) \geq 1,$ then $\card(H_0) \in R_0$
and $\card(H_1) \in R_1$.
\item[(iii)]
If $A$ is an $\mathcal{R}$-basis of order $\mathcal{H}$,
then $A_N = A \cap [0,N]$ is a 
finite $\mathcal{R}$-basis of order $\mathcal{H}$
for every $N \geq 0.$
\item[(iv)]
If $A \neq \{0\}$ is a finite $\mathcal{R}$-basis 
of order $\mathcal{H}$, then $F' = F \setminus \{\max(F)\}$ 
is also a finite $\mathcal{R}$-basis of order $\mathcal{H}$.
\eenum
\et

\pf
To prove~(i) and~(ii), we observe that if $r_A(0,H_0) \geq 1$, then $0 \in A.$
If $r_A(1,H_1) \geq 1$, then $1 \in A.$
Since, for every $h \geq 1,$ both 0 and 1 have unique representations 
as sums of exactly $h$ nonnegative integers,
it follows that if $r_A(0,H_0)\in R_0,$ then
\[
r_{A}(0,H_0) = \sum_{h_0 \in H_0} r_{A}(0,h_0) = \sum_{h_0 \in H_0} 1 = \card(H_0) \in R_0.
\]
Similarly, if $r_A(1,H_1)\in R_1,$ then
\[
r_{A}(1,H_1) = \card(H_1) \in R_1.
\]

The statements~(iii) and~(iv) follow immediately from the definition of a finite basis.
\eop

\section{K\" onig's lemma}
The principal tool in this paper is K\" onig's lemma on the existence
of infinite paths in trees.  For completeness, we include a short proof below.
 
A {\em graph} $G$ consists of a nonempty set $\{v\}$, 
whose elements are called {\em vertices}, 
and a set $\{e\}$, whose elements are called {\em edges}.
Each edge is a set $e = \{v,v'\}$,
where $v$ and $v'$ are vertices and $v \neq v'.$ 
Thus, we are considering only graphs without loops or multiple edges.

We use the following terminology.
The vertices $v$ and $v'$ are called {\em adjacent} if $\{v,v'\}$ is an edge.
The degree of a vertex $v$ is the number of edges $e$ with $v \in e$.
A {\em path} in $G$ from vertex $v$ to vertex $v'$ is a sequence of vertices 
$v_0,v_1, v_2,\ldots, v_n$ such that $v_0 = v$, $v_n = v'$, 
and $v_{i-1}$ is adjacent to $v_{i}$ for all $i = 1,\ldots,n$. 
We define the {\em length} of this path by $n$.
The graph $G$ is {\em connected} if for every two vertices $v$ and $v'$
with $v \neq v'$ there is a path from $v$ to $v'.$
A graph is connected if and only if, for some vertex $v_0$, there is a path
from $v_0$ to $v$ for every vertex $v \neq v_0$.

A {\em simple path} in $G$ is a path whose vertices are pairwise distinct.
A {\em simple circuit} is a sequence of vertices $v_0,v_1, v_2,\ldots, v_n$
such that $n \geq 3$, $\{v_{i-1}, v_{i}\}$ is an edge for $i = 1,\ldots,n$,
$v_i \neq v_j$ for $0 \leq i < j \leq n-1$, and $v_0 = v_n$. 
A graph $G$ has no simple circuits if and only, for every pair of distinct 
vertices $v$ and $v'$, there is at most one simple path from $v$ to $v'$.
An {\em infinite simple path} is an infinite sequence of pairwise distinct vertices 
$v_0, v_1, v_2,\ldots$ such that $v_{i-1}$ is adjacent to $v_{i}$ for all $i \geq 1$. 

A {\em tree} is a connected graph with no simple circuits.
A {\em rooted tree} is a tree with a distinguished vertex, 
called the {\em root} of the tree.
In a rooted tree, for every vertex $v$ different from the root, 
there is a unique simple path in the tree from the root to $v$.

\bt[K\" onig's lemma]
If $T$ is a rooted tree with infinitely many vertices 
such that every vertex has finite degree,
then $T$ contains an infinite simple path beginning at the root.
\et 

\pf
Let $v_0$ be the root of the tree.
We use induction to prove that for every $n$ there is a simple path $v_0, v_1, \ldots, v_n$ 
such that the tree $T$ contains infinitely many vertices $v$ for which the 
unique simple path from the root $v_0$ to $v$ begins with the vertices
$v_0, v_1, \ldots, v_n.$
Since $T$ has infinitely many vertices, the root $v_0$ satisfies this condition.

Let $n\geq 1,$ and assume that we have constructed a simple path 
$v_0, v_1,\ldots, v_{n-1}$ of vertices of the tree $T$ 
with the property that $T$ contains an infinite set $I_{n-1}$ of vertices
such that, for every $v \in I_{n-1}$, the unique simple path 
from $v_0$ to $v$ passes through vertex $v_{n-1}.$
Since the degree of $v_{n-1}$ is finite, the set $F_n$ 
of vertices $v \neq v_{n-2}$ that are adjacent to $v_{n-1}$ is a finite set.
For every vertex $v \in I_{n-1}$, there is a unique simple path in $T$ that 
begins at $v_0$, passes through $v_{n-1}$ and exactly one of the vertices in 
$F_n$, and ends at $v$.  
By the pigeonhole principle, 
since $I_{n-1}$ is infinite, there is a vertex $v_n \in F_n$ 
and an infinite set $I_n \subseteq I_{n-1}$ of vertices 
such that, for every $v \in I_n$, the unique path from $v_0$ to $v$ passes through $v_n.$
This completes the induction.  
The vertices $v_0, v_1, v_2,\ldots$ are pairwise distinct,
and $v_0, v_1, v_2,\ldots$ is an infinite simple path in $T$.
\eop

\section{The generalized representation function problem}

In this section we prove that there exists an
infinite $\mathcal{R}$-basis of order $\mathcal{H}$
if and only if there exist arbitrarily large finite 
$\mathcal{R}$-bases of order $\mathcal{H}$.

\bt  \label{konig:theorem:main}
Let $\mathcal{R} = \{R_n\}_{n=0}^{\infty}$ and 
$\mathcal{H} = \{H_n\}_{n=0}^{\infty}$ 
be sequences of nonempty finite sets of positive integers such that 
\beq             \label{konig:hn}
\lim_{n\rightarrow\infty} \frac{\max(H_n)}{n} = 0.
\eeq
There exists an $\mathcal{R}$-basis $A$ of order $\mathcal{H}$
if and only if for every $N$ there exists a 
finite $\mathcal{R}$-basis $A_N$ of order $\mathcal{H}$
with $\max(A_N) \geq N.$
\et

\pf
If $A$ is a $\mathcal{R}$-basis of order $\mathcal{H}$,
then, by~(\ref{konig:hn}) and Theorem~\ref{konig:theorem:finite}, the set $A$ is infinite, 
hence for every $N$ there is 
an integer $a(N) \in A$ with $a(N) \geq N.$  
By Theorem~\ref{konig:theorem:facts}, the set
$A_N = A \cap [0,a(N)]$ is a finite $\mathcal{R}$-basis of order $\mathcal{H}$
with $\max(A_N) \geq N.$

Conversely, suppose that for every $N$ there exists a 
finite $\mathcal{R}$-basis $A_N$ of order $\mathcal{H}$
with $\max(A_N) \geq N.$  
If $N \geq 1,$ then $0,1\in A_N$ and the sets $\{0\}$ and $\{0,1\}$ are finite 
$\mathcal{R}$-bases of order $\mathcal{H}$.

We construct the graph $T$ whose vertices are the 
finite $\mathcal{R}$-bases of order $\mathcal{H}$.
This graph has infinitely many vertices, since there are 
finite $\mathcal{R}$-bases of order $\mathcal{H}$
with arbitrarily large maximum elements.

Vertices $V$ and $V'$ will be called adjacent in this graph 
if $V' \subseteq V$ and $V \setminus V' = \{\max(V)\}.$
The sets $\{0\}$ and $\{0,1\}$ are adjacent vertices of this graph,
and $\{0,1\}$ is the only vertex adjacent to $\{0\}.$
If $V$ is a vertex and $\card(V) \geq 2,$ then
it follows from Theorem~\ref{konig:theorem:facts}
that $V' = V \setminus \{\max(V)\}$ is a vertex.  Moreover, $V'$ is the 
unique vertex adjacent to $V$ in $T$ such that $\card(V') = \card(V) - 1$.
If $V''$ is adjacent to $V$ and $V'' \neq V'$, then $V = V'' \setminus \{\max(V'')\}$
and $\card(V'') = \card(V)+1.$

We shall prove that $T$ is a rooted tree with root $V_0 = \{0\}$.
Let $V = \{a_0,a_1,\ldots,a_n\}$ be a vertex,
where $0 = a_0 < a_1 < \cdots < a_n.$
For every $k = 0,1,\ldots,n,$ the set $V_k = \{a_0,a_1,\ldots,a_k\}$
is a finite $\mathcal{R}$-basis of order $\mathcal{H}$,
hence is  a vertex of $T$.  Then $V_0 = \{0\}$, $V_n = V,$ and 
$V_0, V_1,\ldots, V_n$ is a simple path in $T$ from the root $V_0$ to $V$.
It follows that the graph $T$ is connected.

Suppose that $n \geq 3$ and $V_0,V_1,\ldots,V_{n-1}, V_n$ is a simple circuit in $T$,
where $V_n = V_0.$  Let $V_{n+1} = V_1.$ 
Since each vertex is a finite set of integers, we can choose $k \in [1,n]$ 
such that $V_k$ is a vertex in the circuit of maximum cardinality.
Vertices $V_{k-1}$ and $V_{k+1}$ are adjacent to $V_k,$
hence $\card(V_k)-\card(V_{k-1})=\pm1$ and $\card(V_k)-\card(V_{k+1})=\pm1$ .
The maximality of $\card(V_k)$ implies that  
$\card(V_k) - \card(V_{k-1}) = \card(V_k) - \card(V_{k+1}) = 1$,
and so $V_{k-1} = V_k \setminus\{\max(V_k)\} = V_{k+1},$ 
which is impossible.
Therefore, $T$ contains no simple circuit, and so $T$ is a tree.

To apply K\" onig's lemma, we must prove 
that every vertex of this tree has finite degree.
The only vertex adjacent to the root $\{0\}$ is $\{0,1\}$,
hence $\{0\}$ has finite degree.  
Let $V \neq \{0\}$ be a vertex of $T$.
Then $1 \in V$ and so $\max(V) \geq 1.$
Suppose that the $V$ is adjacent to infinitely many vertices $V'$.
The only subset of $V$ that is a vertex adjacent to $V$ is $V \setminus \{\max(V)\}.$
Every other vertex $V'$ adjacent to $V$ is a superset of $V$
of the form $V' = V \cup \{\max(V')\}.$
For each such $V'$, the integer $n = \max(V')-1$ must be an element of the sumset  
$h_nV$ for some $h_n \in H_n$, and so
\[
n \leq h_n\max(V) \leq \max(H_n)\max(V).
\]
Since $n \in h_nV$ for infinitely many integers $n$, 
it follows that
\[
\limsup_{n\rightarrow\infty}\frac{\max(H_n)}{n} \geq \frac{1}{\max(V)} > 0,
\]
which contradicts~(\ref{konig:hn}).
Thus, every vertex of the infinite tree $T$ has finite degree.
By K\" onig's lemma, the tree must contain an infinite simple path 
$\{0\} = V_0, V_1, V_2, \ldots.$
For each nonnegative integer $n$, let $a_n = \max(V_n).$  
Then $V_n = \{a_0,a_1,\ldots,a_n\}$ for all $n = 0,1,2,\ldots.$
Let
\[
A = \{a_n\}_{n=0}^{\infty} = \bigcup_{n=0}^{\infty} V_n.
\]
Since $a_n \geq n,$ it follows that 
\[
r_{A}(n,H_n)  = \sum_{h_n \in H_n} r_A(n,h_n) = \sum_{h_n \in H_n} r_{V_n}(n,h_n) 
= r_{V_n}(n,H_n) \in R_n,
\]
and so $A$ is a $\mathcal{R}$-basis of order $\mathcal{H}$.
This completes the proof.
\eop

\bt
Let $h \geq 2$ and let $f$ be a function such that $f(n)$ 
is a positive integer for every nonnegative integer $n.$
There exists a basis $A$ of order $h$ with representation function $r_A(n,h) = f(n)$
if and only if for every $N$ there exists a finite set $A_N$
of nonnegative integers with $\max(A_N) \geq N$ and $r_A(n,h) = f(n)$
for all $n = 0, 1,\ldots,\max(A_N).$
\et

\pf
This follows immediately from Theorem~\ref{konig:theorem:main} with $H_n = \{h\}$
and $R_n = \{ f(n)\}$ for all nonnegative integers $n$.
\eop

Applying Theorem~\ref{konig:theorem:main} 
to the classical Erd\H os-Tur\' an conjecture, we obtain
the following result of Dowd~\cite[Theorem 2.1]{dowd88}.

\bt
Let $c \geq 1$ and $h \geq 2$.  There exists a basis $A$ of order $h$ such that  
\[
r_{A}(n,h) \leq c \qquad\text{for all $n \geq 0$}
\]
if and only if, for every $N$,
there exists a finite set $A_N$  of nonnegative integers 
such that $\max(A_N) \geq N$ 
and
\[
1 \leq r_{A_N}(n,h) \leq c
\]
for all $n = 0,1,\ldots,\max(A_N).$
\et

\pf
This follows immediately from Theorem~\ref{konig:theorem:main}
with $H_n = \{h\}$ and $R_n = [1,c]$ for all $n \geq 0$.
\eop

\section{Ordered representation functions}
There are other important representation functions 
in additive number theory.
For example, for any set $A$ of integers, 
the {\em ordered representation function}  
$r'_A(n,h)$ counts the number of $h$-tuples $(a_1,\ldots,a_h) \in A^h$ such that
$a_1 + \cdots + a_h = n$. 
Nathanson~\cite{nath78c} proved the following uniqueness theorem 
for ordered representation functions:  For any function $f:\N_0 \rightarrow\N_0$
and for any positive integer $h$, there exists at most one set $A$ of 
nonnegative integers such that $r'_A(n,h) = f(n)$ for all $n \in \N_0$.
He also showed that uniqueness does not hold if ordered representation 
functions only eventually coincide, and he described all pairs 
of sets $A$ and $B$ of nonnegative integers 
such that $r'_A(n,2) = r'_B(n,2)$ for all sufficiently large integers $n$.

We can define {\em basis of order $\mathcal{H}$} and 
{\em $\mathcal{R}$-basis of order $\mathcal{H}$}
in terms of the ordered representation function.
Theorem~\ref{konig:theorem:main} is also true for ordered representation functions.


\begin{thebibliography}{1}

\bibitem{dowd88}
M.~Dowd, \emph{Questions related to the {Erd\H os-Tur\' an} conjecture}, SIAM
  J. Discrete Math. \textbf{1} (1988), 142--150.

\bibitem{erdo-tura41}
P.~Erd\H{o}s and P.~Tur\'an, \emph{On a problem of {S}idon in additive number
  theory and some related questions}, J. London Math. Soc. \textbf{16} (1941),
  212--215.

\bibitem{grek-hadd-helo-pihk03}
G.~Grekos, L.~Haddad, C.~Helou, and J.~Pihko, \emph{On the {Erd\H os-Tur\' an}
  conjecture}, Preprint, 2002.

\bibitem{nath78c}
M.~B. Nathanson, \emph{Representation functions of sequences in additive number
  theory}, Proc. Amer. Math. Soc. \textbf{72} (1978), 16--20.

\bibitem{nath03c}
M.~B. Nathanson, \emph{Every function is the representation function of an additive
  basis for the integers}, Preprint, 2003.

\end{thebibliography}
\end{document}